\documentclass[12pt]{article}

\usepackage{amsfonts,amsmath,amssymb}

\begin{document}

\begin{center}
  {\large \bf On the local factor of the zeta 
function of quadratic orders}
  \vspace{2.5ex}\\
  {\sc Masanobu Kaneko}\\
\end{center}

\begin{abstract}
  We prove by an elementary method the Riemann hypothesis for the
  local Euler factor of the zeta function of quadratic orders.
\end{abstract}

\footnotetext{\!\!\!\!{\it 2000 Mathematics Subject Classification}: 
Primary 11R42; Secondary 11M38.

}

Let $K$ be a quadratic number field, ${\cal O}_K$ its ring of
integers, and ${\cal O}_f$ for each integer $f\ge1$ the order of
conductor $f$ in ${\cal O}_K={\cal O}_1$.
As shown in \cite{K} and \cite{Z}, the zeta function $\zeta_{{\cal
    O}_f}(s)$  of ${\cal O}_f$, which is defined by
$$\zeta_{{\cal O}_f}(s)=\sum_{\frak a}\frac1{N(\frak a)^s},$$
where the
sum extends over all {\it proper} ideals $\frak a$ of ${\cal O}_f$
with norm $N(\frak a)$, has the following form:
$$
\zeta_{{\cal O}_f}(s)=\zeta_K(s)\cdot\prod_{p\vert f}
\frac{(1-p^{-s})(1-\chi(p)p^{-s})-p^{n_p-1-2n_ps}
  (1-p^{1-s})(\chi(p)-p^{1-s})}{(1-p^{1-2s})}.
$$
Here, $\zeta_K(s)$ is the Dedekind zeta function of $K$, the
product is over the prime factors $p$ of the conductor $f$, with $n_p$
being the exact power of $p$ in $f$, and $\chi$ is the  Dirichlet
character corresponding to the extension $K/{\mathbf Q}$.  \vspace{2ex}

The main purpose of this short note is to provide proofs of the
following properties, including the ``Riemann hypothesis,'' for the local factor
$$
\varepsilon_{f,\,p}(s):=\frac{(1-p^{-s})(1-\chi(p)p^{-s})-p^{n_p-1-2n_ps}
  (1-p^{1-s})(\chi(p)-p^{1-s})}{(1-p^{1-2s})}
$$
of $\zeta_{{\cal O}_f}(s)/\zeta_K(s)$: 
\vspace{2ex}

\noindent {\bf Theorem.} 1) \begin{em}
  The function $\varepsilon_{f,\,p}(s)$ is a polynomial of degree $2n_p$ in
  $p^{-s}$ and  satisfies the functional equation
$$\varepsilon_{f,\,p}(1-s)=p^{-n_p(1-2s)}\varepsilon_{f,\,p}(s).$$ \end{em}

2) \begin{em} All the zeros of $\varepsilon_{f,\,p}(s)$ lie on the
  line $\mbox{Re} (s)=1/2$.
\end{em}
\vspace{2ex}

\noindent {\it Proof.} ~ Setting $u=p^{-s}$ and $n=n_p$, we rewrite 
$\varepsilon_{f,\,p}(s)$ as the function $P_n(u)$, given as follows:
\begin{eqnarray*} 
P_n(u)&=&\frac{(1-u)(1-\chi(p)u)-p^{n-1}u^{2n}(1-pu)(\chi(p)-pu)}
{1-pu^2}\\
&=& \frac{1-p^{n+1}u^{2(n+1)}-(1+\chi(p))u(1-p^nu^{2n})+\chi(p)u^2
(1-p^{n-1}u^{2(n-1)})}{1-pu^2}. 
\end{eqnarray*}
The numerator of this expression vanishes if we set $u=\pm 1/\sqrt{p}$
and hence is
divisible by the denominator $1-pu^2$. Thus $P_n(u)$ is indeed a
polynomial of degree $2n$. By direct division, we find
$$P_n(u)=1-(1+\chi(p))u+\cdots
-p^{n-1}(1+\chi(p))u^{2n-1}+p^nu^{2n}.$$
The functional equation can
be verified straightforwardly. This ends the proof of assertion 1).

To prove the Riemann hypothesis 2), put $s=1/2+it/\log p$. Then we have
$u=p^{-1/2}e^{-it}$ and
$$P_n(u)=\frac{1-e^{-2(n+1)it}-(1+\chi(p))p^{-1/2}e^{-it}(1-e^{-2nit})
  +\chi(p)p^{-1}e^{-2it}(1-e^{-2(n-1)it})}{1-e^{-2it}}.$$
Then, using the relation
$$\frac{1-e^{-2mit}}{1-e^{-2it}}=\frac{e^{-mit}}{e^{-it}} \frac{\sin
  mt}{\sin t},$$
we obtain
$$P_n(u)=\frac{e^{-nit}}{p}\left(p\,\frac{\sin (n+1)t}{\sin t}-
  \sqrt{p}\,(1+\chi(p))\,\frac{\sin nt}{\sin t}-\chi(p)\,\frac{\sin
    (n-1)t}{\sin t}\right).$$
We have to show that if the  right-hand
side of this is zero then $t$ is real.  Recall that, for any integer
$m\ge0$, the quotient $\sin (m+1)t/\sin t$ is a polynomial of degree
$m$ in $x=\cos t$, which is referred to as the Chebyshev polynomial of
the second kind, denoted by $U_m(x)$.  The first several of these are
as follows:
$$U_0(x)=1,\quad U_1(x)=2x,\quad U_2(x)=4x^2-1,\quad U_3(x)=8x^3-4x.$$
Note that the function $U_m(x)$ can also be defined for $m<0$; in
particular, we have $U_{-1}(x)=0$ and $U_{-2}(x)=-1$.  Using the
$U_m(x)$, the proof is reduced to showing that all the roots of the
polynomial (of degree $n$)
$$Q_n(x):=p\,U_n(x)-\sqrt{p}\,(1+\chi(p))\,U_{n-1}(x)+\chi(p)\,U_{n-2}(x)
\quad (n\ge1)$$
are in the real interval $[-1,1]$. 

Because of the recurrence of the Chebyshev polynomials
$$U_n(x)=2xU_{n-1}(x)-U_{n-2}(x),$$
the polynomials $Q_n(x)$ satisfy the same
recurrence:
$$Q_n(x)=2x Q_{n-1}(x)-Q_{n-2}(x)\quad (n\ge 2),$$
with
$Q_0(x)=p-\chi(p)$. We show that the $n$ roots of $Q_n(x)$ are all in
the interval $(-1,1)$ by making use of the theorem of Sturm
(cf. \cite[\S92]{W}), utilizing the fact that $Q_n(x),Q_{n-1}(x),\ldots, Q_0(x)$ 
forms a ``Strum sequence''. Because $U_n(1)=n+1$ and $U_n(-1)=(-1)^n(n+1)$, we have
\begin{eqnarray*}
  Q_n(1)&=&p(n+1)-\sqrt{p}(1+\chi(p))n+\chi(p)(n-1)\\
&=& (\sqrt{p}-1)(\sqrt{p}-\chi(p))n+p-\chi(p)>0
\end{eqnarray*}
and
\begin{eqnarray*}
  Q_n(-1)&=&p(-1)^n(n+1)-\sqrt{p}(1+\chi(p))(-1)^{n-1}n+\chi(p)
(-1)^{n-2}(n-1)\\
&=&(-1)^n\left\{p(n+1)+\sqrt{p}(1+\chi(p))n+\chi(p)(n-1)\right\}\\
&=&(-1)^n\left\{(\sqrt{p}+1)(\sqrt{p}+\chi(p))n+p-\chi(p)\right\}.
\end{eqnarray*}
The sign of the last expression is $(-1)^n$, and hence the number of
``variations,'' as defined in \cite[\S92]{W}, is $n$.  Then, noting the theorem
of Sturm, we conclude that $Q_n(x)$ has $n$ roots in the
interval $(-1,1)$. (Note that, since the  degree of $Q_n(x)$ is $n$,
condition 4 of \cite[\S92]{W} need not to be verified.) \hfill\hbox{\rule{5pt}{8pt}}\medskip

\noindent
{\bf Remarks and questions.} 1) It is amusing that the properties
stated in the above theorem are precisely those enjoyed by the
congruence zeta function (or, rather, its essential part) of a curve
of genus $n=n_p$ over the prime field $\mathbf{F}_p$. This naturally
leads us to wonder  if $\varepsilon_{f,\,p}(s)$ admits a cohomological
(or any other ``nice'') interpretation and if the above theorem can be
proved ``conceptually'' using such an interpretation.

2) The theorem proved here shows in particular that the quotient
$\zeta_{{\cal O}_f}(s)/ \zeta_K(s)$ is entire and that the Riemann
hypothesis holds for $\zeta_{{\cal O}_f}(s)$ only if it holds for
$\zeta_K(s)$. It is known that for a Galois  extension $k'/k$ of
number fields, the quotient of the Dedekind zeta functions
$\zeta_{k'}(s)/\zeta_k(s)$ is entire. Thus the zeta function of the
{\it over} field is divisible by that of the base field.  On the
contrary, in the case considered above, the zeta function of the
{\it subring}  ${\cal O}_f$ is divisible by that of the over ring. What
is the reason for this?

3) The generating function of the polynomials $P_n(x)$ takes the simple form
$$F(u,X):=1+\sum_{n=1}^\infty P_n(u)X^n=\frac{(1-uX)(1-\chi(p)uX)}
{(1-X)(1-pu^2X)},$$
and the functional equation for $P_n(u)$, which is written as
$$p^nu^{2n}P_n(\frac1{pu})=P_n(u),$$ is encoded as
$$F(\frac1{pu},pu^2X)=F(u,X).$$

4)  For another zeta function
$$\zeta_{{\cal O}_f}^{\ast}(s)=\sum_{\frak a}\frac1{N(\frak a)^s},$$
where $\frak a$ runs over {\it all} (not necessarily proper) ideals of
${\cal O}_f$, the Euler product is ({\it cf.} \cite{Z})
$$\zeta_K(s)\cdot\prod_{p\vert f}\frac{1-p^{(n_p+1)(1-2s)}-\chi(p)
  p^{-s}(1-p^{n_p(1-2s)})}{1-p^{1-2s}}.$$
It can be shown similarly
that the local factor $$\frac{1-p^{(n_p+1)(1-2s)}-\chi(p)
  p^{-s}(1-p^{n_p(1-2s)})}{1-p^{1-2s}}$$
possesses the same properties as in the theorem.

5) It would be nice to have a generalization of our theorem to the
zeta functions of orders of number fields of higher degree.

\noindent {\bf Acknowledgment.} The author is grateful to Christopher
Deninger for his interest in the present work, without which this paper
would not have been written.



\vspace{3ex}
\noindent Graduate School of Mathematics, Kyushu University 33, \\
Fukuoka, 812-8581 JAPAN\\
{\tt mkaneko@math.kyushu-u.ac.jp}\\


\begin{thebibliography}{9}
  
\bibitem{K} M.~Kaneko : A generalization of the  Chowla-Selberg
  formula and the zeta functions of quadratic orders,  {\em
    Proc. Japan Acad.}, {\bf 66(A)-7} (1990), 201--203.
  
\bibitem{W} H.~Weber : {\em Lehrbuch der Algebra}, Vol. 1, Chelsea,
  New York.
  
\bibitem{Z} D.~Zagier : Modular forms whose Fourier  coefficients
  involve zeta functions of quadratic fields,  in {\it Modular
    functions of one variable VI}, Lect. Notes in Math., no. 627,
  Springer-Verlag, (1977) 105--169.

\end{thebibliography}
\end{document}